\documentclass[a4paper,reqno,12pt]{amsart}
\usepackage{rotating,pdflscape,color,amssymb,hyperref,subfigure,psfrag,amsmath,eufrak,bbm}
\author[Florent Benaych-Georges]{Florent Benaych-Georges}\address{Florent Benaych-Georges, LPMA,  UPMC Univ Paris 6, Case courier 188, 4, Place Jussieu, 75252 Paris Cedex 05, France, and CMAP, \'Ecole Polytechnique, route de Saclay, 91128 Palaiseau Cedex, France.} \email{florent.benaych@upmc.fr}
\urladdr{http://www.cmapx.polytechnique.fr/\~\/benaych/}

\author[Raj Rao Nadakuditi]{Raj Rao Nadakuditi}\address{Raj Rao Nadakuditi, Department of Electrical Engineering and Computer Science, University of Michigan, 1301 Beal Avenue, Ann Arbor, MI 48109. USA.}
\email{rajnrao@eecs.umich.edu}
\urladdr{http://www.eecs.umich.edu/\~\/rajnrao/}

\title[Low rank perturbations of large random matrices]{The singular values and  vectors of low rank perturbations of large rectangular random matrices}

\thanks{F.B.G's work was partially supported by the \emph{Agence Nationale de la Recherche} grant ANR-08-BLAN-0311-03. R.R.N's work was supported by an Office of Naval Research Young Investigator Award N00014-11-1-0660.}

\date{\today}
\newcommand{\lto}{\longrightarrow}

\newcommand{\ti}{\widetilde}
\newcommand{\bet}{\beta}
\newcommand{\al}{\alpha}
\newcommand{\tta}{\theta}
\newcommand{\Tta}{\Theta}
\newcommand{\si}{\sigma}
\newcommand{\udl}{\underline}

\newcommand{\ovl}{\overline}
\newcommand{\bbm}{\begin{bmatrix}}
\newcommand{\ebm}{\end{bmatrix}}
\newcommand{\bes}{\begin{equation*}}
\newcommand{\ees}{\end{equation*}}
\newcommand{\be}{\begin{equation}}
\newcommand{\ee}{\end{equation}}
\newcommand{\beqy}{\begin{eqnarray}}
\newcommand{\eeqy}{\end{eqnarray}}
\newcommand{\beq}{\begin{eqnarray*}}
\newcommand{\eeq}{\end{eqnarray*}}

\newcommand{\lan}{\langle}
\newcommand{\ran}{\rangle}

\newcommand{\diag}{\operatorname{diag}}
\newcommand{\Diag}{\operatorname{diag}}

\newcommand{\convas}{\overset{\textrm{a.s.}}{\longrightarrow}}
\newcommand{\convind}{\overset{\mathcal{D}}{\longrightarrow}}

\newcommand{\Pro}{\mathbb{P}}

\newcommand{\Span}{\operatorname{Span}}

\newcommand{\Tr}{\operatorname{Tr}}

\newcommand{\ninf}{\underset{n\to\infty}{\longrightarrow}}

\newcommand{\one}{\mathbbm{1}}

\newcommand{\R}{\mathbb{R}}
\newcommand{\C}{\mathbb{C}}

\newcommand{\K}{\mathbb{K}}

\newcommand{\ud}{\mathrm{d}}

\newcommand{\pro}{probability }
\newcommand{\ie}{\textit{i.e.}}
\newcommand{\f}{\frac}
\newcommand{\ff}{\frac{1}}
\newcommand{\lf}{\left}
\newcommand{\ri}{\right}

\newcommand{\st}{such that }

\newcommand{\vfi}{\varphi}
\newcommand{\ste}{\textrm{ s.t. }}

\newcommand{\eps}{\varepsilon}

\newcommand{\bxp}{\boxplus}

\newcommand{\wtX}{\widetilde{X}_n}

\newcommand{\wtu}{\widetilde{u}}
\newcommand{\wtv}{\widetilde{v}}
\newcommand{\wts}{\widetilde{\si}_1}

\newcommand{\bck}{\backslash}

\newtheorem{Th}{Theorem}[section]

\newtheorem{hyp}[Th]{Assumption}

\newtheorem{propo}[Th]{Proposition}

\newtheorem{lem}[Th]{Lemma}

\newtheorem{rmq}[Th]{Remark}

\newenvironment{pr}{\noindent {\it Proof. }}{\hfill$\square$}

\long\def\symbolfootnote[#1]#2{\begingroup
\def\thefootnote{\fnsymbol{footnote}}\footnote[#1]{#2}\endgroup}

\begin{document}
\maketitle
\begin{abstract}In this paper,  we consider the singular values  and singular vectors of finite,
low rank perturbations of large rectangular random matrices. Specifically, we prove
almost sure convergence of the extreme singular values and appropriate
projections of the corresponding singular vectors of the perturbed matrix.

As in the prequel, where we considered the eigenvalues of Hermitian matrices, the non-random limiting  value is shown to depend explicitly on the limiting singular value distribution of the unperturbed matrix via an integral transform that linearizes rectangular additive convolution in free probability theory. The asymptotic position  of the extreme singular values of the perturbed matrix differs from that of the original matrix if and only if the singular values of the perturbing matrix are above a certain critical threshold which depends on this same aforementioned integral transform.

We examine the consequence of this singular value phase transition on the associated left and right singular eigenvectors and discuss the   fluctuations around these non-random limits.
\end{abstract}




\section{Introduction}
In many applications, the  $n \times m$ signal-plus-noise data or measurement matrix formed by stacking the $m$ samples or measurements of $n \times 1$ observation vectors alongside each other can be modeled as:
\begin{equation}\label{eq:model}
\widetilde{X} = \sum_{i=1}^{r} \sigma_{i} u_{i} v_{i}^{*} + X,
\end{equation}
where $u_{i}$ and $v_{i}$ are left and right ``signal'' column vectors, $\sigma_{i}$ are the associated ``signal'' values and $X$ is the noise-only matrix of random noises. This model is ubiquitous in signal processing \cite{vantrees02a,scharf91a}, statistics \cite{Muirhead82,anderson84a,jolliffeprincipal} and machine learning \cite{kannan2009spectral} and is known under various guises as a signal subspace model \cite{schmidt86a}, a latent variable statistical model \cite{jordan1998learning}, or a probabilistic PCA model \cite{tipping1999probabilistic}.

Relative to this model, a common application-driven objective is to estimate the signal subspaces Span$\{u_{1}, \ldots, u_{r}\}$ and Span$\{v_{1}, \ldots, v_{r}\}$ that contain signal energy. This is accomplished by computing the \emph{singular value decomposition} (in brief \emph{SVD}) of $\widetilde{X}$ and extracting the $r$ largest singular values and the associated singular vectors of $\widetilde{X}$ - these are referred to as the $r$ \textit{principal components} \cite{rao1964use} and the   Eckart-Young-Mirsky theorem states that they provide the best rank-$r$ approximation of the matrix $\widetilde{X}$ for any unitarily invariant norm \cite{eckart1936approximation,mirsky1960symmetric}. This theoretical justification combined with the fact that these vectors can be efficiently computed using now-standard numerical algorithms for the SVD \cite{golub1965calculating} has led to the ubiquity of the SVD in applications such as array processing \cite{vantrees02a}, genomics \cite{alter2000singular}, wireless communications \cite{edfors2002ofdm}, information retrieval \cite{furnas1988information} to list a few \cite{klema2002singular}.

In this paper, motivated by emerging high-dimensional statistical applications \cite{johnstone2009statistical}, we place ourselves in the setting where $n$ and $m$ are large and the SVD of $\widetilde{X}$ is used to form estimates  of $\{\sigma_{i}\}$, $\{u_{i}\}_{i=1}^{r}$ and $\{v_{i}\}_{i=1}^{r}$. We provide a   characterization of the relationship between the estimated extreme singular values of $\widetilde{X}$ and the true ``signal'' singular values $\sigma_{i}$ (and also the angle between the estimated and true singular vectors). 

In the limit of large matrices, the extreme singular values only depend on integral transforms of the distribution of the singular values of the noise-only matrix $X$ in \eqref{eq:model} and exhibit a phase transition about a critical value: this is a new occurrence of the so-called \emph{BBP phase transition}, named after the authors of the seminal paper \cite{bbp05}.  The critical value   also depends on the   aforementioned integral transforms which arise from rectangular free probability theory \cite{b09, bg07c}.
We also characterize the fluctuations of the singular values about these asymptotic limit. The results obtained are precise in the large matrix limit and, akin to our results in \cite{benaych-rao.09}, go beyond answers that might be obtained using matrix perturbation theory \cite{ss90}.

Our results are in a certain sense very general (in terms of possible distributions for the noise model $X$) and recover as a special case results found in the literature for the eigenvalues \cite{bbp05,bs06}  and eigenvectors \cite{hr07,p07,n08}  of $\widetilde{X}\widetilde{X}^{*}$ in the setting where $X$ in (\ref{eq:model}) is Gaussian. For the Gaussian setting we provide new results for the right singular vectors.
Such results had already been proved in the particular case where $X$ is a Gaussian matrix, but our approach brings to light a general principle, which can be applied beyond the Gaussian case. Roughly speaking, this principle says that for $X$ a   $n\times p$ matrix (with $n,p\gg 1$), if one adds an independent  small rank perturbation  $ \sum_{i=1}^{r} \sigma_{i} u_{i} v_{i}^{*} $ to $X$, then the extreme singular values will   move to  positions which are approximately the solutions $z$ of the equations $$\ff{n}\Tr \f{z}{z^2I-XX^*}\times \ff{p}\Tr \f{z}{z^2I-X^*X}=\ff{\tta_i^2},\qquad\textrm{($1\le i\le r$)}.$$

In the case where these equations have no solutions (which means that the $\tta_i$'s are below a certain threshold), then the extreme  singular values of $X$ will not move significantly. We also provide similar results for the associated left and right singular vectors and give limit theorems for the fluctuations.  These expressions provide the basis for the parameter estimation algorithm developed by Hachem et al in \cite{hachemloubatonmnv11}.

 The  papers  \cite{benaych-rao.09,bg-maida-guionnetTCL} were devoted to the analogue problem for the eigenvalues of finite rank perturbations of Hermitian matrices. We follow the strategy developed in these papers for our proofs: we derive master equation representations that implicitly encode the relationship between the singular values and singular vectors of $X$ and $\widetilde{X}$ and use concentration results   to obtain the stated analytical expressions. Of course, because of these similarities  in the proofs, we chose to focus, in the present paper, in what differs from  \cite{benaych-rao.09,bg-maida-guionnetTCL}. 

 At a certain level, our proof also present analogies with the ones of other papers devoted to other occurrences of the BBP phase transition, such as \cite{sandrinePTRF06,FP07,CDF09,CDF09b,CDFF10}. We mention that the approach of the paper   \cite{bg-maida-guionnetPGD} could also be used to consider large deviations of the extreme singular values of $\widetilde{X}$.

This paper is organized as follows. We state our main results in Section \ref{assumptions.1272010} and provide some examples in Section \ref{sec:examples}. The proofs are provided in Sections \ref{sec:proof 1}-\ref{sec:proof fluct} with some technical details relegated to the appendix in Section \ref{sec:appendix}.

\section{Main results}\label{assumptions.1272010}
\subsection{Definitions and hypotheses}\label{def_and_hyp}
Let $X_n$ be a $n \times m$ real or complex random matrix. Throughout this paper we assume that  $n \leq m $ so that we may simplify the exposition of the proofs. We may do so without loss of generality because in the setting where $n > m$, the expressions derived will hold for $X_n^{*}$.  Let the $n \leq m$ singular values of $X_n$ be $\si_{1} \geq \si_{2} \geq \ldots \geq \si_n$. Let $\mu_{X_n}$ be the empirical singular value distribution, \ie, the probability measure defined as
$$ \mu_{X_n} = \dfrac{1}{n} \sum_{i=1}^n \delta_{\si_i}.$$
Let $m$ depend on $n$ -- we denote this dependence explicitly by $m_n$  which we will sometimes omit for brevity by substituting $m$ for $m_n$. Assume that as $n \longrightarrow \infty$,   $n/m_n \longrightarrow c \in [0,1]$. In the following, we shall need some of the following hypotheses.

\begin{hyp}\label{hypspec}The probability measure $\mu_{X_n}$ converges almost surely weakly  to a non-random compactly supported probability measure $\mu_X$.\end{hyp}

 Examples of random matrices satisfying this   hypothesis can be found in e.g. \cite{bai-silver-95, capcas.iumj, bai-silver-book, b09, agz09, pan}. Note however that the question of isolated extreme singular values is not addressed in papers like \cite{bai-silver-95,pan} (where moreover the perturbation has a non bounded rank).

\begin{hyp}\label{hypspeca}  Let   $a$ be infimum of the support of $\mu_X$. The smallest singular value of $X_n$ converges almost surely to  $a$.\end{hyp}

\begin{hyp}\label{hypspecb} Let   $b$ be supremum of the support of $\mu_X$. The largest singular value of $X_n$ converges almost surely to  $b$.\end{hyp}

Examples of random matrices satisfying the above   hypotheses can be found in e.g. \cite{capdonmart, bai-silver-book,  agz09, pan}.

In this problem, we shall consider the extreme singular values and the associated singular vectors of $\wtX$, which is the random $n\times m$ matrix:
$$\wtX=X_n+P_n,$$
where $P_n$ is defined as described below.

For a given $r \ge 1$, let $\tta_1\ge\cdots\ge \tta_r >0$ be deterministic non-zero real numbers, chosen independently of $n$.    For every $n$, let   $G_u^{(n)}, G_v^{(n)}$ be two independent matrices with sizes respectively $n\times r$ and $m\times r$, with i.i.d. entries distributed according to a fixed \pro measure $\nu$ on $\K=\R$ or $\C$. We introduce the column vectors  $u_1, \ldots, u_r\in \K^{n\times 1}$ and $ v_1,�\ldots , v_r\in \K^{m\times 1}$  obtained from $G_{u}^{(n)}$ and $G_{v}^{(n)}$ by either:
\begin{enumerate}
\item Setting $u_{i}$ and $v_{i}$ to equal the $i$-th column of $\ff{\sqrt{n}}G_u^{(n)}$ and $\ff{\sqrt{m}}G^{(n)}_v$ respectively or,
\item Setting $u_{i} $ and $v_{i}$ to equal to the vectors obtained from a Gram-Schmidt (or QR factorization) of $G_u^{(n)}$ and $G^{(n)}_v$ respectively.
\end{enumerate}
We shall refer to the model (1) as   the {\it i.i.d. model} and to the model (2) as the {\it orthonormalized model}.
 With the $u_i$'s and $v_i$'s constructed as above, we define the random perturbing matrix $P_n\in \K^{n\times m}$ as:
 $$P_n = \sum_{i=1}^{r} \theta_{i} u_{i}  v_{i}^{*}.$$ In the orthonormalized model, the $\tta_i$'s are the non zero singular values of $P_n$ and the $u_i$'s and the $v_i$'s are the left and right associated singular vectors.

We make the following hypothesis on the law $\nu$ of the entries of  $G_u^{(n)}$ and $ G_v^{(n)}$ (see \cite[Sect. 2.3.2]{agz09} for the definition of log-Sobolev inequalities).
\begin{hyp} \label{hyponG}
 The \pro measure $\nu$ has  mean zero, variance one and that satisfies a log-Sobolev inequality.\end{hyp}

\begin{rmq}{\rm
We also note if $\nu$ is the standard real or complex Gaussian distribution, then the singular vectors produced using the orthonormalized model will have uniform distribution on the set of $r$ orthogonal random vectors.}
\end{rmq}

\begin{rmq}{\rm
If $X_n$ is random but has a bi-unitarily invariant distribution and $P_n$ is non-random with rank $r$, then we are in same setting as the orthonormalized model for the results that follow. More generally, our idea in defining both of our models (the i.i.d. one and the orthonormalized one) was to show that if $P_n$ is chosen independently from $X_n$ in a somehow ``isotropic way'' (i.e. via a distribution which is not faraway from being invariant by the action of the orthogonal group by conjugation), then a BBP phase transition occurs, which is governed by a certain integral transform of the limit empirical singular values distribution of $X_n$, namely $\mu_X$.}
\end{rmq}

\begin{rmq}{\rm
We note that there is small albeit non-zero probability that $r$ i.i.d. copies of a random vector are not linearly independent. Consequently, there is a small albeit non-zero probability that the $r$ vectors obtained as in (2) via the Gram-Schmidt orthogonalization may not be well defined. However, in the limit of large matrices, this process produces well-defined vectors with overwhelming probability (indeed, by Proposition \ref{concentration_U.9610}, the determinant of the associated $r\times r$ Gram matrix tends to one). This is implicitly assumed in what follows.}
\end{rmq}

\subsection{Notation} Throughout this paper, for $f$ a function and $d \in \R$, we set $$f(d^+):=\lim_{z\downarrow d}f(z)\,;\qquad f(d^-):=\lim_{z\uparrow d}f(z),$$
we also let $\convas$ denote almost sure convergence. The (ordered) singular values of an $n\times m$ Hermitian matrix $M$ will be denoted by $\si_1(M)\ge\cdots \ge \si_n(M)$. Lastly, for a subspace $F$ of a Euclidian space $E$ and a unit vector $x\in E$, we denote the norm of the orthogonal projection of $x$ onto $F$ by $\lan x, F\ran$.

\subsection{Largest singular values and singular vectors phase transition}
In Theorems \ref{140709.main.rectangular}, \ref{180709.13h39.rectangular} and \ref{deloc.sing.vect.250709.rectangular}, we suppose Assumptions \ref{hypspec}, \ref{hypspecb} and \ref{hyponG} to hold.

We define $\ovl{\tta}$, the threshold of the phase transition, by the formula $$\ovl{\tta}:=(D_{\mu_X}(b^+))^{-1/2},$$ with the convention that $(+\infty)^{-1/2}=0$, and where
 $D_{\mu_X}$, the $D$-transform of $\mu_X$ is the function, depending on $c$, defined by
$$ D_{\mu_X}(z) :=\left[\int \frac{z}{z^2-t^{2}} \ud \mu_X(t)\right] \times \left[c\int \frac{z}{z^2-t^{2}} \ud \mu_X(t)+\frac{1-c}{z}\right] \qquad \textrm{for } z>b,$$
In the theorems below,  $D_{\mu_X}^{-1}(\cdot)$ will denote its functional inverse on $[b,+\infty)$.

\begin{Th}[Largest singular value phase transition]\label{140709.main.rectangular}
The $r$ largest singular values of the $n \times m$ perturbed matrix
$\wtX$
exhibit the following behavior as $n,m_n \to \infty$ and $n/m_n \to c$. We have that for each fixed $1\leq i \leq r$,
$$\si_i(\wtX) \convas \begin{cases} D_{\mu_X}^{-1}(1/\theta_i^2)&\textrm{ if $\theta_i>\ovl{\tta}$,}\\ \\
b &\textrm{ otherwise.}\end{cases}$$
Moreover, for each fixed $i>r$, we have that $\si_i(\wtX)\convas b$.
\end{Th}

\begin{Th}[Norm of projection of largest singular vectors]\label{180709.13h39.rectangular}Consider indices ${i_0}\in \{1, \ldots, r\}$   \st $\theta_{i_0}>\ovl{\tta}  $
. For each $n$, define $\widetilde{\si}_{i_0}=\si_{i_0}(\wtX)$
and let $\widetilde{u}$ and $\widetilde{v}$ be   left and right unit singular vectors of $\wtX$ associated with the singular value $\widetilde{\si}_{i_0}$. Then we have, as $n \longrightarrow \infty$,

\flushleft a)
\be\label{250709.09h13}
|\langle \widetilde{u}, \Span\{u_i\ste \tta_i=\tta_{i_0}\}\rangle|^{2}  \convas {\f{-2\vfi_{\mu_X}(\rho)}{\theta_{i_0}^2 D'_{\mu_X}(\rho) }},
\ee
\flushleft b)
\be\label{250709.09h13.300709}
|\langle \widetilde{v}, \Span\{v_i\ste \tta_i=\tta_{i_0}\}  \rangle|^{2}  \convas  {\f{-2\vfi_{\widetilde{\mu}_X}(\rho)}{\theta_{i_0}^2 D_{\mu_X}'(\rho) }},
\ee
 where $\rho=D_{\mu_X}^{-1}(1/\theta_{i_0}^2)$ is the limit of $\widetilde{\si}_{i_0}$ and $\widetilde{\mu}_X=c\mu_X+(1-c)\delta_0$ and for any \pro measure $\mu$,
\be\label{9710.19h07}\vfi_\mu(z):=\int \f{z}{z^2-t^2}\ud \mu(t).\ee
\flushleft c) Furthermore, in the same asymptotic limit, we have
$$ |\langle \widetilde{u}, \Span\{u_i\ste \tta_i\ne \tta_{i_0}\} \rangle|^{2}    \convas 0, \quad \textrm{ and } \quad |\langle \widetilde{v},  \Span\{v_i\ste \tta_i\ne \tta_{i_0}\} \rangle|^{2} \convas 0 ,$$
and
$$\lan \vfi_{{\mu}_X}(\rho){P_n}\widetilde{v}-\widetilde{u} \;,\;  \Span\{u_i\ste \tta_i=\tta_{i_0}\} \rangle \convas 0.$$
\end{Th}

\begin{Th}[Largest singular vector phase transition]\label{deloc.sing.vect.250709.rectangular}When $r=1$, let the  sole singular value of $P_n$  be denoted by $\theta$. Suppose that
\be\label{25709.10h08}
\theta\le \ovl{\tta} \qquad\textrm{and} \qquad\vfi_{\mu_X}'(b^+)=-\infty.
\ee
For each $n$, let $\widetilde{u}$ and $\widetilde{v}$ denote, respectively,  left and right unit singular vectors of $\wtX$  associated with its largest singular value. Then
$$ \langle \widetilde{u},   \ker (\theta^2 I_n-P_nP_n^*)\rangle \convas 0, \quad \textrm{ and } \quad\langle \widetilde{v},   \ker (\theta^2I_m-P_n^*P_n)\rangle \convas 0,$$
as $n \longrightarrow \infty$.
\end{Th}

The following proposition allows to assert that in many classical matrix models,  the threshold $\ovl{\tta}$ of the  above phase transitions  is positive.   The proof relies on a straightforward computation which we omit.

\begin{propo}[Edge density decay condition for phase transition]\label{square root add}
Assume that the limiting singular distribution $\mu_X$ has a density $f_{\mu_X}$ with a power decay at $b$, i.e., that, as $t\to b$ with $t<b$,  $f_{\mu_X}(t) \sim M (b-t)^\alpha$  for some exponent $\alpha>-1$ and  some constant $M$. Then:
$$ \ovl{\tta}=(D_{\mu_X}(b^{+}))^{-1/2} >0\iff \alpha>0  \qquad\textrm{ and } \qquad \vfi'_{\mu_X}(b^{+})=-\infty \iff \al\le 1,$$
so that the phase transitions in Theorems \ref{140709.main.rectangular} and \ref{deloc.sing.vect.250709.rectangular} manifest for $\al=1/2$.
\end{propo}

\begin{rmq}[Necessity of singular value repulsion for the singular vector phase transition]\label{210709.23h19}{\rm
Under additional hypotheses on the manner in which the empirical singular distribution of $X_n \convas \mu_X$ as $n \longrightarrow \infty$,  Theorem \ref{deloc.sing.vect.250709.rectangular} can be generalized to any singular value with limit    $b$ \st $D_{\mu_X}'(\rho)$ is infinite. The specific hypothesis has to do with requiring the spacings between the singular values of $X_n$  to be more ``random matrix like'' and exhibit repulsion instead of being ``independent sample like'' with possible clumping. We plan to develop this line of inquiry in a separate paper.}
\end{rmq}

\subsection{Smallest singular values and vectors for square matrices}\label{smallest_sing_sect}We now consider the phase transition exhibited by the smallest singular values and vectors. We restrict ourselves to the setting where $\wtX$ is a square matrix; this restriction is necessary because the non-monotonicity of the function  $D_{\mu_X}$ on $[0, a)$ when   $c =\lim n/m <1$, poses some technical difficulties that do not arise in the square setting.
Moreover, in Theorems \ref{smallest_eigenvalue_rect}, \ref{sing_vect_square} and \ref{smallest_sing_vect_square}, we suppose Assumptions \ref{hypspec}, \ref{hypspeca} and \ref{hyponG} to hold.

We define $\udl{\tta}$, the threshold of the phase transition, by the formula $$\udl{\tta}:=(\vfi_{\mu_X}(a^-))^{-1},$$ with the convention that $(+\infty)^{-1}=0$, and where $\vfi_{\mu_X}(z)=\int\f{z}{z^2-t^2}\ud\mu(t)$, as  in Equation \eqref{9710.19h07}.
In the theorems below,  $\vfi_{\mu_X}^{-1}(\cdot)$ will denote its functional inverse of the function $\vfi_{\mu_X}(\cdot)$  on $(0, a)$.

\begin{Th}[Smallest singular value phase transition for square matrices]\label{smallest_eigenvalue_rect}
When $a>0$ and $m =n $,  the $r$ smallest singular values of $\wtX$ exhibit the following behavior. We have that for each  fixed $1 \leq i \leq r$,
$$\si_{n+1-i}(\wtX) \convas
\begin{cases}
\vfi_{\mu_X}^{-1}(1/\theta_i)&\textrm{ if  $\theta_i>\udl{\tta}$,}\\ \\ a &\textrm{ otherwise.}
\end{cases}$$
Moreover, for each fixed $i>r$, we have that $\si_{n+1-i}(\wtX)\convas a$.
\end{Th}

\begin{Th}[Norm of projection of smallest singular vector for square matrices]\label{sing_vect_square}
Consider indices ${i_0}\in \{1, \ldots, r\}$   \st $\theta_{i_0}>\udl{\tta}$. For each $n$, define $\widetilde{\si}_{i_0}=\si_{n+1-i_0}(\wtX)$
and let $\widetilde{u}$ and $\widetilde{v}$ be   left and right unit singular vectors of $\wtX$ associated with the singular value $\widetilde{\si}_{i_0}$. Then we have, as $n \longrightarrow \infty$,

\flushleft a)
\be\label{250709.09h13square}
|\langle \widetilde{u}, \Span\{u_i\ste \tta_i=\tta_{i_0}\}\rangle|^{2}  \convas \f{-1}{\vfi'_{\mu_X}(\rho)},
\ee
\flushleft b)
\be\label{250709.09h13.300709square}
|\langle \widetilde{v}, \Span\{v_i\ste \tta_i=\tta_{i_0}\}  \rangle|^{2}  \convas  \f{-1}{\vfi'_{\mu_X}(\rho)},
\ee
\flushleft c) Furthermore, in the same asymptotic limit, we have
$$ |\langle \widetilde{u}, \Span\{u_i\ste \tta_i\ne \tta_{i_0}\} \rangle|^{2}    \convas 0, \quad \textrm{ and } \quad |\langle \widetilde{v},  \Span\{v_i\ste \tta_i\ne \tta_{i_0}\} \rangle|^{2} \convas 0 ,$$
and
$$\lan \vfi_{{\mu}_X}(\rho){P_n}\widetilde{v}-\widetilde{u} \;,\;  \Span\{u_i\ste \tta_i=\tta_{i_0}\} \rangle \convas 0, $$\end{Th}

\begin{Th}[Smallest singular vector phase transition]\label{smallest_sing_vect_square}
When $r = 1$ and $m = n$, let the smallest singular value of $\wtX$ be denoted by $\widetilde{\si}_{n}$ with $\widetilde{u}$ and $\widetilde{v}$ representing  associated left and right unit singular vectors respectively. Suppose that
\bes\label{25709.10h082}
a>0, \qquad \theta \leq \udl{\tta} \qquad\textrm{and} \qquad\vfi_{\mu_X}'(a^-)=-\infty.
\ees
Then
$$ \langle \widetilde{u},   \ker (\theta^2 I_n-P_nP_n^*)\rangle \convas 0, \quad \textrm{ and } \quad\langle \widetilde{v},   \ker (\theta^2I_m-P_n^*P_n)\rangle \convas 0,$$
as $n \longrightarrow \infty$.
\end{Th}

The analogue of Remark \ref{210709.23h19} also applies here.

\subsection{The $D$-transform in free probability theory}

The {\it $C$-transform with ratio $c$} of a \pro measure $\mu$ on $\R_+$,  defined as:
\be\label{rel_D_H}C_{\mu}(z)=U\lf( {z}(D_\mu^{-1}(z))^2-1\ri),\ee
where the function $U$, defined as:
$$
U(z)=  \begin{cases} \f{-c-1+\lf[(c+1)^2+4cz\ri]^{1/2}}{2c} & \textrm{ when } c>0 ,\\ z  & \textrm{ when }  $c = 0$, \end{cases}
$$
is the analogue of the logarithm of the Fourier transform for the {\it rectangular free  convolution with ratio $c$} (see \cite{b09R=C, bg.sph.int} for an introduction to the theory of rectangular free convolution) in the sense described next.

Let $A_n$ and $B_n$ be independent $n \times m$ rectangular random matrices that are invariant, in law, by conjugation by any orthogonal (or unitary) matrix. Suppose that, as $n, m \to \infty$ with $n/m \to c$, the empirical singular values distributions $\mu_{A_n}$ and $\mu_{B_n}$ of $A_n$ and $B_n$ satisfy $\mu_{A_{n}} \longrightarrow \mu_{A}$ and $\mu_{B_{n}} \longrightarrow \mu_{B}$. Then by \cite{b09},  the empirical singular values distribution $\mu_{A_n + B_n}$ of $A_n+B_n$ satisfies $\mu_{A_n + B_n} \longrightarrow \mu_{A} \bxp_{c} \mu_{B}$, where $ \mu_{A} \bxp_{c} \mu_{B}$ is  a \pro measure which
 can be characterized in terms of the $C$-transform as
$$C_{\mu_A\bxp\mu_B}(z)= C_{\mu_A}(z)+C_{\mu_B}(z).$$
The coefficients of the series expansion  of $U(z)$ are the   {\it rectangular free cumulants with ratio $c$} of $\mu$ (see \cite{bg07c} for an introduction to the  rectangular free cumulants). The connection between free  rectangular additive convolution and $D_\mu^{-1}$ (via the $C$-transform) and the appearance of  $D_\mu^{-1}$ in Theorem \ref{140709.main.rectangular} could be of independent interest to free probabilists: the emergence of this transform in the study of isolated singular values completes the picture of \cite{benaych-rao.09}, where the transforms linearizing additive and multiplicative free convolutions already appeared in similar contexts.

\subsection{Fluctuations of the largest singular value}

Assume that the empirical singular value distribution of $X_n$ converges to $\mu_X$ faster than $1/\sqrt n.$ More precisely,
\begin{hyp}\label{hyp_fluctu_rect} We have  $$\f{n}{m_n}=c+o(\ff{\sqrt{n}}),$$ $r=1$, $\tta:=\tta_1>\ovl{\tta} $ and
   $$\ff{n}\Tr(\rho^2 I_{n} -X_nX_n^*)^{-1}=\int\f{1}{\rho^2-t^2}\ud{\mu_X}(t) +o(\ff{\sqrt{n}})$$ for $\rho=D_{\mu_X}^{-1}(1/\tta^2)$ the limit of $\si_1(\wtX)$.
\end{hyp}

We also make the following hypothesis on the law $\nu$ (note that it doesn't contains the fact that $\nu$ is symmetric). In fact, wouldn't it hold, we would still have a limit theorem on the fluctuations of the largest singular value, like in Theorem 3.4 of \cite{bg-maida-guionnetTCL}, but we chose not to develop this case.
\begin{hyp} \label{hyponGbis}
 If  $\nu$ is entirely supported by the real line, $\int x^4 \ud \nu(x)=3$. If  $\nu$ is not entirely supported by the real line, the real and imaginary parts of a $\nu$-distributed random variables are
 independent and identically distributed with $\int |z|^4 \ud \nu(z)=2$. \end{hyp}

 Note that we do not ask  $\nu$ to be symmetric and make no hypothesis about its third moment. The reason is that the main ingredient of the following theorem is Theorem 6.4 of \cite{bg-maida-guionnetTCL} (or Theorem 7.1 of \cite{bai-yao-TCL}), where no hypothesis of symmetry or about the third moment is done.

\begin{Th}\label{sinfluoutbulk} Suppose Assumptions \ref{hypspec}, \ref{hypspecb},  \ref{hyponG}, \ref{hyp_fluctu_rect} and  \ref{hyponGbis} to hold.
Let $\widetilde{\si}_1$ denote the largest singular value of $\wtX$. Then as $n \longrightarrow \infty$,
$$ n^{1/2} \left(\widetilde{\si}_1 - \rho \right) \convind \mathcal{N}(0,s^{2}),$$
where $\rho = D_{\mu_X}^{-1}(c, 1/\tta^2)$ and
$$s^2=\begin{cases}\dfrac{f^{2}}{2\beta}&\textrm{for the i.i.d. model,}\\ \\
\dfrac{f^{2} - 2}{2\beta}&\textrm{for the orthonormalized model,}
\end{cases}$$
with $\beta = 1 $ (or $2$) when $X$ is real (or complex) and
$$f^{2} := \f{\int \f{\ud \mu_X(t)}{(\rho^2-t^2)^2}}{\lf(\int \f{\ud \mu_X(t)}{\rho^2-t^2}\ri)^2}+\f{\int \f{\ud \widetilde{\mu}_X(t)}{(\rho^2-t^2)^2}}{\lf(\int \f{\ud \widetilde{\mu}_X(t)}{\rho^2-t^2}\ri)^2}+2\f{\int \f{t^2\ud \mu_X(t)}{(\rho^2-t^2)^2}}{ \int \f{\rho\ud \mu_X(t)}{\rho^2-t^2}\int \f{\rho\ud \widetilde{\mu}_X(t)}{\rho^2-t^2}},$$
with $\widetilde{\mu}_X =  c\mu_X+(1-c)\delta_0$.
\end{Th}

\subsection{Fluctuations of the smallest singular value of square matrices} When $m_n = n$ so that $c = 1$, assume that:
\begin{hyp}\label{H2prime} For all $n$, $m_n=n$,  $r=1$, $\tta:=\tta_1>\udl{\tta}$ and
   $$ \ff{n}\Tr(\rho^2 I_n-X_nX_n^*)^{-1}=\int\f{1}{\rho^2-t^2}\ud{\mu_X}(t)+o(\ff{\sqrt{n}})$$ for $\rho:=\vfi_{\mu_X}^{-1}(1/\tta)$ the limit of the smallest singular value of $\widetilde{X}_n$.
\end{hyp}

\begin{Th}\label{sinfluoutbulksquare2}Suppose Assumptions \ref{hypspec}, \ref{hypspeca},  \ref{hyponG}, \ref{H2prime} and  \ref{hyponGbis} to hold. Let $\widetilde{\si}_n$ denote the smallest singular value of  $\widetilde{X}_n$.
Then as $n \longrightarrow \infty$
$$n^{1/2} \left(\widetilde{\si}_n-\rho \right) \convind \mathcal{N}(0,s^{2}),$$
where
$$s^{2} = \begin{cases} \dfrac{f^2}{2\beta} & \textrm{ for the i.i.d. model} \\ \\
\dfrac{f^{2}-2}{2\beta} & \textrm{ for the orthonormalized model}
\end{cases}$$
with $\beta = 1 $ (or $2$) when $X$ is real (or complex) and
$f^{2} := 2 \theta^{2} \int \f{\kappa^2+t^2}{(\kappa^2-t^2)^2}\ud \mu_X(t)$.
\end{Th}

\section{Examples}\label{sec:examples}
\subsection{Gaussian rectangular random matrices with non-zero mean}
Let $X_n$ be an $n\times m$ real (or complex) matrix with independent, zero mean, normally distributed entries with variance $1/m$. It is known \cite{mp67, bai-silver-book} that, as $n,m \longrightarrow \infty$ with $n/m\to c \in (0,1]$,  the spectral measure of the singular values of $X_n$  converges to the 
 distribution with density
\[
\ud \mu_{X}(x)= \f{\sqrt{4c-(x^2-1-c)^2}}{\pi c x}\one_{(a,b)}(x) dx, \]
where $a=1-\sqrt{c}$ and $b=1+\sqrt{c}$ are the end points of the support of $\mu_X$.  It is known \cite{bai-silver-book}  that the extreme eigenvalues converge to the bounds of this support.

Associated with this singular measure, we have, by an application of the result in \cite[Sect. 4.1]{fbg05.inf.div} and  Equation \eqref{rel_D_H},
\beq  &  D_{\mu_{X}}^{-1}(z)=\sqrt{\f{(z+1)(cz+1)}{z}}\, ,&\\ &D_{\mu_{X}}(z)=\f{z^2-(c+1)-\sqrt{(z^2-(c+1))^2-4c}}{2c}\, ,\quad D_{\mu_{X}}(b^+)=\ff{\sqrt{c}}\, .&\eeq

Thus for any $n\times m$ deterministic matrix $P_{n}$ with $r$ non-zero singular values $\theta_1\geq \cdots\geq \theta_r$ ($r$ independent of $n,m$), for any fixed $i\geq 1$, by Theorem \ref{140709.main.rectangular}, we have
 \be
 \si_i(X_n+P_n) \convas
 \begin{cases}
 \sqrt{\frac{(1+\theta_i^2)(c+\theta_i^2)}{\theta_{i}^{2}}} & \textrm{if $i\le r$ and } \theta_{i} > c^{1/4} \\
 1+\sqrt{c} & \textrm{otherwise}.
 \end{cases}
 \ee
 as $n\longrightarrow \infty$. As far as the i.i.d. model is concerned, this formula allows us to recover some of the results of \cite{bbp05}.

 Now, let us turn our attention to the singular vectors. In the setting where $r=1$,  let $P_n = \theta u v^{*}$. Then, by Theorems \ref{180709.13h39.rectangular} and \ref{deloc.sing.vect.250709.rectangular}, we have
\be\label{290909.20h07} \left|\lan \widetilde{u}, u\ran \right|^2\convas  \begin{cases}1-\dfrac{c(1+\theta^2)}{\theta^2(\theta^2+c)}
&\textrm{if $\theta \geq {c}^{1/4}$,}
\\
0&\textrm{otherwise.}
\end{cases}\ee

The phase transitions for the eigenvectors of $\wtX^*\wtX$ or for the pairs of singular vectors of $\wtX$ can be similarly computed to yield the expression:

\be\label{290909.20h07.v} \left|\lan \widetilde{v}, v \ran \right|^2\convas  \begin{cases}1 - \dfrac{(c +\theta^{2})}{\theta^{2}(\theta^{2} +1 )}
&\textrm{if $\theta \geq {c}^{1/4}$,}
\\
0&\textrm{otherwise.}
\end{cases}\ee

\subsection{Square Haar unitary matrices}
Let $X_n$ be Haar distributed unitary (or orthogonal) random matrix. All of its singular values are equal to one, so that it has limiting spectral measure
$$
\mu_{X}(x) = \delta_1,
$$
with $a = b = 1$ being the end points of the support of $\mu_X$.

Associated with this spectral measure, we have (of course, $c=1$)
$$D_{\mu_{X}}(z)=\f{z^2}{(z^2-1)^2}\qquad \textrm{for } z\geq 0, \; z\neq 1,$$ thus for all $\theta>0$, $$D_{\mu_{X}}^{-1}( 1/\theta^2)=\begin{cases}\f{\theta+\sqrt{\theta^2+4}}{2}&\textrm{if the inverse is computed on $(1, +\infty)$,}\\
 \f{-\theta+\sqrt{\theta^2+4}}{2}&\textrm{if the inverse is computed on $(0,1)$.}\end{cases}$$

Thus for any $n \times n$, rank $r$ perturbing matrix $P_n$ with  $r$ non-zero singular values $\theta_1\geq  \cdots\geq \theta_r$ where neither $r$, nor the  $\theta_i$'s  depend on $n$, for any fixed $i=1, \ldots, r$, by Theorem \ref{140709.main.rectangular} we have
\[
\si_i(X_n+P_n) \convas \f{\theta_i+\sqrt{\theta_i^2+4}}{2}
\quad \textrm{ and }\quad
\si_{n+1-i}(X_n+P_n) \convas\f{-\theta_i+\sqrt{\theta_i^2+4}}{2}\]
while for any fixed    $i\geq r+1$, both $\si_i(X_n+P_n)$ and  $\si_{n+1-i}(X_n+P_n)$ $\convas 1$.

\section{Proof of Theorems \ref{140709.main.rectangular} and \ref{smallest_eigenvalue_rect}} \label{sec:proof 1}
The proofs of both theorems are quite similar. As a consequence, we only prove Theorem \ref{140709.main.rectangular}.

The sequence of steps described below yields the desired proof (which is very close to the one of Theorem 2.1 of \cite{benaych-rao.09}):
\begin{enumerate}
\item{The first, rather trivial, step in the proof of Theorem \ref{140709.main.rectangular} is to use Weyl's interlacing inequalities to prove that any fixed-rank singular value of $\wtX$ which does not tend to a limit $>b$  tends   to   $b$.}
\item{Then, we utilize Lemma \ref{agl_lin_8710} below  to express the extreme singular values of $\wtX$ as the $z$'s \st   a certain random $2r\times 2r$  matrix $M_n(z)$ is singular.}
\item{We then exploit convergence properties of certain analytical functions (derived in the appendix) to prove that almost surely,  $M_n(z)$ converges  to a certain  deterministic matrix $M(z)$, uniformly in $z$.}
\item{We then invoke a continuity lemma (see Lemma \ref{roots_continuity_9610}  in the appendix) to claim that almost surely, the $z$'s \st $M_n(z)$ is singular ({\it i.e.} the extreme singular values of $\wtX$) converge to the $z$'s \st $M(z)$ is singular.}
\item{We conclude  the proof by noting that, for our setting, the $z$'s \st $M(z)$ is singular are precisely the $z$'s \st for some $i\in \{1,\ldots, r\}$, $D_{\mu_X}(z)=\ff{\tta^2_i}$.  Part (ii) of Lemma \ref{roots_continuity_9610} , about the rank of $M_n(z)$, will be useful to assert that when the $\tta_i$'s are pairwise distinct, the multiplicities of the isolated singular values  are all equal to one.}
\end{enumerate}

Firstly, up to a conditioning by the $\si$-algebra generated by the $X_n$'s, one can suppose them to be deterministic  and all the randomness supported by the perturbing matrix $P_n$.

Secondly, by \cite[Th. 3.1.2]{hj91}, one has, for all $i\ge 1$, $$\si_{i+r}(X_n)\le \si_i(\widetilde{X}_n)\le \si_{i-r}(X_n)$$with the convention $\si_j(X_n)=+\infty$ for $i\le 0$ and $0$ for $i>n$.  By the same proof  as in \cite[Sect. 6.2.1]{benaych-rao.09}, it follows that for all $i\ge 1$ fixed, \be\label{moka.1111}\liminf \si_i(\wtX)\ge b\ee and that for all  fixed $i>r$, \be\label{moka.1111.2}\si_i(\wtX)\ninf b\ee
 (we insist here on the fact that $i$ has to be fixed, i.e. not to depend on $n$: of course, for $i=n/2$, \eqref{moka.1111.2} is not true anymore in general).

Our approach is based on the following lemma, which reduces the problem to the study of $2r\times 2r$ random matrices. Recall that the constants $r$, $\tta_1, \ldots, \tta_r$,  and the random column vectors (which depend on $n$, even though  this dependence does not appear in the notation) $u_1, \ldots, v_r$, $v_1, \ldots, v_r$ have been introduced in Section \ref{def_and_hyp} and that the perturbing matrix $P_n$ is given by  $$P_n=\sum_{i_1}^r \tta_iu_iv_i^*.$$ Recall also that the singular values of $X_n$ are denoted by $\si_1\ge \cdots \ge \si_n$.  Let us define the matrices $$\Theta=\diag(\tta_1,\ldots, \tta_r)\in \R^{r\times r}, \quad U_n=\bbm u_1  \cdots u_r \ebm\in \K^{n\times r},\quad V_m=\bbm v_1 \cdots v_r \ebm\in \K^{m\times r}.$$

\begin{lem}\label{agl_lin_8710}
The positive singular values of $\wtX$ which are not singular values of $X_n$ are the  $z\notin \{\si_1, \ldots, \si_n\}$ \st the $2r\times 2r$ matrix$$M_n(z):=\bbm U_n^*(z^2 I_n-X_nX_n^*)^{-1}U_n& U_n^*(z^2 I_n-X_nX_n^*)^{-1}X_n V_m \\ V_m^*X_n^*(z^2 I_n-X_nX_n^*)^{-1}U_n & V_m^*(z^2 I_m-X_n^*X_n)^{-1}V_m \ebm-\bbm 0&\Tta^{-1}\\ \Tta^{-1}& 0\ebm $$ is not invertible.
\end{lem}

For the sake of completeness, we provide a proof, even though several related results can be found in the literature (see e.g. \cite{agg88,bg-maida-guionnetPGD}).

\begin{pr}Firstly, \cite[Th. 7.3.7]{hj85} states that the non-zero singular values of $\wtX$ are the positive eigenvalues of $\bbm 0&\wtX\\ \wtX^*&0\ebm$. Secondly, for any $z>0$ which is   not a singular value of $X_n$,
by \cite[Lem. 6.1]{bg-maida-guionnetTCL}, $$\det\lf(zI_n-\bbm 0&\wtX\\ \wtX^*&0\ebm\ri)= \det\lf(zI_n-\bbm 0&X_n\\ X_n^*&0\ebm\ri)^{-1}\times \prod_{i=1}^r\tta_i^2  \times\det M_n(z),$$ which allows to conclude, since by hypothesis, $\det\lf(z I_{n+m}-\bbm 0&X_n\\ X_n^*&0\ebm\ri)^{-1}\ne 0$.
\end{pr}

Note that  by Assumption \ref{hypspec},  \beq \ff{n}\Tr \f{z}{z^2 I_n-X_nX_n^*} &\ninf&  \int \f{z}{z^2-t^2}\ud\mu_X(t), \\ \ff{m}\Tr \f{z}{z^2 I_m -X_n^*X_n} &\ninf&  \int \f{z}{z^2-t^2}\ud\ti{\mu}_X(t)\qquad (\ti{\mu}_X=c\mu_X+(1-c)\delta_0),\eeq uniformly on
 any subset of $\{z\in \C\ste \Re(z)>b+\eta \}$, $\eta>0$.
It follows, by a direct application of Ascoli's Theorem and Proposition \ref{concentration_U.9610},  that almost surely, we have the following convergence (which is uniform in $z$) \beq\label{8710.14h50} U_n^*\f{z}{z^2 I_n-X_nX_n^*}U_n&\ninf&  \lf(\int \f{z}{z^2-t^2}\ud\mu_X(t)\ri)\cdot I_r, \\
\label{8710.14h51}V_m^*  \f{z}{z^2 I_m-X_n^*X_n}V_m&\ninf&  \lf( \int \f{z}{z^2-t^2}\ud\ti{\mu_X}(t)\ri)\cdot I_r.\eeq In the same way, almost surely \beq\label{8710.14h50.1} U_n^*(z^2 I_n-X_nX_n^*)^{-1}X_n V_m\ninf 0&\textrm{ and }&
 V_m^*X_n^*(z^2 I_n-X_nX_n^*)^{-1}U_n\ninf 0 .\eeq
It follows that almost surely,
\be\label{M_ntoM9610}M_n(z)\ninf M(z):=\bbm \vfi_{\mu_X}(z)I_r&0\\ 0&\vfi_{\widetilde{\mu}_X}(z)I_r\ebm-\bbm 0& \Tta^{-1}\\
\Tta^{-1}&0\ebm, \ee where $\vfi_{\mu_X}$ and $\vfi_{\widetilde{\mu}_X}$ are the functions defined in the statement of Theorem \ref{180709.13h39.rectangular}.

Now, note that once \eqref{moka.1111} has been established,  our result only concerns the number of singular values
of $\wtX$  in $[b+\eta, +\infty)$ (for any $\eta>0$), hence can be proved via Lemma \ref{roots_continuity_9610}. Indeed,  by Hypothesis \ref{hypspecb}, for $n$ large enough, $X_n$ has no singular value $>b+\eta$, thus numbers $>b+\eta$ cannot be in the same time singular values of $X_n$ and $\wtX$.

In the case where the $\tta_i$'s are pairwise distinct, Lemma \ref{roots_continuity_9610} allows to conclude the proof of Theorem \ref{140709.main.rectangular}. Indeed, Lemma \ref{roots_continuity_9610} says that exactly  as much singular values  of $\wtX$ as predicted by the theorem  have limits $>b$ and that their limits are exactly the ones predicted by the Theorem. The part of the theorem devoted to singular values tending to $b$ can then be deduced  from \eqref{moka.1111} and \eqref{moka.1111.2}.

In the case where the $\tta_i$'s are not pairwise distinct, an approximation approach allows to conclude (proceed for example as in Section 6.2.3 of \cite{benaych-rao.09}, using \cite[Cor. 7.3.8 (b)]{hj85} instead of
 \cite[Cor. 6.3.8]{hj85}).

\section{Proof of Theorems \ref{180709.13h39.rectangular} and \ref{sing_vect_square}} \label{sec:proof 2}

The proofs of both theorems are quite similar. As a consequence, we only prove Theorem  \ref{180709.13h39.rectangular}.

As above, up to a conditioning by the $\si$-algebra generated by the $X_n$'s, one can suppose them to be deterministic and all the randomness supported by the perturbing matrix $P_n$.

Firstly, by the Law of Large Numbers, even in the {\it i.i.d. model}, the $u_i$'s and the $v_i$'s are almost surely asymptotically orthonormalized. More specifically, for all $i\ne j$, $$ \lan u_i,u_j
\ran\ninf \one_{i=j}$$ (the same being true for the $v_i$'s). As a consequence, it is enough to prove that

 a')
\be\label{250709.09h13_iid}
\sum_{i\ste \tta_i=\tta_{i_0}}|\langle \wtu,  u_i \rangle |^{2}  \convas {\f{-2\vfi_{\mu_X}(\rho)}{\theta_{i_0}^2D'_{\mu_X}(\rho) }},
\ee

 b')
\be\label{250709.09h13.300709_iid}
\sum_{i\ste \tta_i=\tta_{i_0}}|\langle \wtv,    v_i \rangle |^{2}  \convas  {\f{-2\vfi_{\widetilde{\mu}_X}(\rho)}{\theta_{i_0}^2 D'_{\mu_X}(\rho)  }},
\ee

c')
\be\label{250709.09h13.300709_iid2}\sum_{i\ste \tta_i\ne \tta_{i_0}}| \langle \wtu,  u_i \rangle |^2  +|\langle  \wtv,    v_i   \rangle|^2 \convas 0,\ee

 d')
\be\label{250709.09h13.300709_iid3}\sum_{i\ste \tta_i=\tta_{i_0}}|  \vfi_{{\mu_X}}(\rho)\tta_{i_0}\lan     \wtv,v_i\ran -\lan  \wtu,u_i   \rangle |^2 \convas 0, \ee

Again, the proof is based on a lemma which reduces the problem to the study of the kernel of a random $2r\times 2r$ matrix. The matrices $\Theta$, $U_n$ and $V_m$ are the ones introduced before Lemma \ref{agl_lin_8710}.

\begin{lem}\label{9710.Bob_Dylan_forever}
Let $z$ be a   singular value of $\wtX$ which is not a singular value of $X_n$ and let $(u,v)$ be a corresponding singular pair of unit vectors. Then the column vector \bes\label{sa_eyed_lady}\bbm \Theta{V_m}^*v\\  \Theta{U_n}^*u\ebm\ees belongs to the kernel of the $2r\times 2r$ matrix $M_n(z)$ introduced in Lemma \ref{agl_lin_8710}. Moreover,
we have \beqy \nonumber v^*P_n^*\f{z^2}{(z^2 I_n-X_nX_n^*)^2}P_nv+u^*P_n\f{X_n^*X_n}{(z^2 I_m-X_n^*X_n)^2}P_n^*u  &&\\
+v^*P_n^*\f{z}{(z^2 I_n-X_nX_n^*)^2}X_nP_n^*u
+u^*P_nX_n^*\f{z}{(z^2 I_n-X_nX_n^*)^2}P_nv&=&1\label{11h4412071978}.\eeqy
\end{lem}

\begin{pr}The  first part of the lemma is easy to verify with the formula $X_n^*f(X_nX_n^*)=f(X_n^*X_n)X_n^*$ for any function $f$ defined on $[0,+\infty)$. For the second part, use   the formulas $$\wtX\wtX^* u=z^2u\quad\textrm{ and }\quad X_n^*u=zv-P_n^*u,$$ to establish $
u=(z^2 I_n-X_nX_n^*)^{-1}(zP_nv+X_nP_n^*u),
$
and then use the fact that $u^*u=1$.\end{pr}

Let us consider $z_n$, $(\wtu,\wtv)$ as in the statement of Theorem \ref{180709.13h39.rectangular}. Note firstly that for $n$ large enough, $z_n>\si_1(X_n)$, hence Lemma \ref{9710.Bob_Dylan_forever} can be applied, and the vector \be\label{grizzly_bear_11710} \bbm \Theta {V_m}^*\wtv\\  \Theta {U_n}^*\wtu\ebm=\bbm\tta_1 \lan v_1,\wtv\ran,\ldots,\tta_r \lan v_r,\wtv\ran,\tta_1 \lan u_1,\wtu\ran,\ldots,\tta_r \lan u_r,\wtu\ran \ebm^T\ee belongs to $\ker M_n(z_n)$. As explained in the proof of Theorem  \ref{140709.main.rectangular}, the random matrix-valued function $M_n(\cdot)$ converges almost surely uniformly to the matrix-valued function $M(\cdot)$ introduced in Equation  \eqref{M_ntoM9610}. Hence $M_n(z_n)$ converges almost surely   to $M(\rho)$, and  it follows that the orthogonal projection on $(\ker M(\rho))^\perp$ of the vector of \eqref{grizzly_bear_11710}   tends almost surely to zero.

Let us now compute this projection. For $x,y$ column vectors of $\K^r$, \beq M(\rho)\bbm x\\ y \ebm=0&\iff&\forall i,\quad y_i=\tta_i\vfi_{\mu_X}(\rho)x_i\quad\textrm{ and }\quad x_i=\tta_i\vfi_{\widetilde{\mu}_X}(\rho)y_i\\
&\iff& \forall i, \quad\begin{cases}x_i=y_i=0&\textrm{ if }\quad\tta_i^2\vfi_{\mu_X}(\rho)\vfi_{\widetilde{\mu}_X}(\rho)\ne 1,\\
 y_i=\tta_i\vfi_{ {\mu}_X}(\rho)x_i& \textrm{ if } \quad\tta_i^2\vfi_{\mu_X}(\rho)\vfi_{\widetilde{\mu}_X}(\rho)= 1.\end{cases}
 \eeq
 Note that $\rho$ is precisely defined by the relation $\tta_{i_0}^2\vfi_{\mu_X}(\rho)\vfi_{ \widetilde{\mu}_X}(\rho)= 1$. Hence with   $\bet:=-\tta_{i_0}\vfi_{\widetilde{\mu}_X}(\rho)$, we have,
$$
\ker M(\rho)= \{\bbm x\\ y\ebm\in \K^{r+r}\ste  \forall i, \quad x_i=y_i=0\textrm{ if $\tta_i\neq\tta_{i_0}$ and }  y_i=-\bet x_i \textrm{ if $\tta_i=\tta_{i_0}$}\},
$$ hence
$$(\ker M(\rho))^\perp= \{\bbm x\\ y\ebm\in \K^{r+r}\ste  \forall i, \quad x_i=\bet y_i\quad \textrm{ if $\tta_i=\tta_{i_0}$}\} $$ and the orthogonal projection of any vector $\bbm x\\ y\ebm$ on $(\ker M(\rho))^\perp$ is the vector $\bbm x'\\ y'\ebm$ \st for all $i$, $$(x'_i,y'_i)=\begin{cases}
(x_i,y_i)&\textrm{if $\quad \tta_i\ne\tta_{i_0}$,}\\
\f{ \bet x_i+ y_i}{\bet^2+1} ( \bet,1 )&\textrm{if $\quad \tta_i=\tta_{i_0}$.}
\end{cases}$$

Then,  \eqref{250709.09h13.300709_iid2} and \eqref{250709.09h13.300709_iid3} are  direct consequences of the fact that
the projection of the vector of \eqref{grizzly_bear_11710} on $(\ker M(\rho))^\perp$ tends to zero.

Let us now prove \eqref{250709.09h13.300709_iid}.  By \eqref{11h4412071978}, we have \be\label{12710.20h}a_n+b_n+c_n+d_n=1,\ee with \beqy \label{12710.13h1} &&a_n= {\wtv}^*P_n^*\f{z_n^2}{(z_n^2 I_n-X_nX_n^*)^2}P_n{\wtv}=\sum_{i,j=1}^r\tta_i\tta_j\ovl{\lan v_i, {\wtv}\ran}\lan v_j, {\wtv}\ran{u_i}^*\f{z_n^2}{(z_n^2 I_n-X_nX_n^*)^2}u_j  \\
\label{12710.13h2}&&b_n={\wtu}^*P_n\f{X_n^*X_n}{(z_n^2 I_m-X_n^*X_n)^2}P_n^*{\wtu} =\sum_{i,j=1}^r\tta_i\tta_j\ovl{\lan u_i, {\wtu}\ran}\lan u_j, {\wtu}\ran{v_i}^*\f{X_n^*X_n}{(z_n^2 I_m-X_n^*X_n)^2}v_j \\
 \nonumber&&c_n={\wtv}^*P_n^*\f{z_n}{(z_n^2 I_n-X_nX_n^*)^2}X_nP_n^*{\wtu}=\sum_{i,j=1}^r\tta_i\tta_j\ovl{\lan v_i, {\wtv}\ran}\lan u_j, {\wtu}\ran{u_i}^*\f{z_n}{(z_n^2 I_n-X_nX_n^*)^2}X_nv_j \\
\nonumber&&d_n={\wtu}^*P_nX_n^*\f{z_n}{(z_n^2 I_n-X_nX_n^*)^2}P_n{\wtv}=\sum_{i,j=1}^r\tta_i\tta_j\ovl{\lan u_i, {\wtu}\ran}\lan v_j, {\wtv}\ran{v_i}^*X_n^*\f{z_n}{(z_n^2 I_n-X_nX_n^*)^2}u_j\eeqy
Since the limit of $z_n$ is out of the support of $\mu_X$, one can apply Proposition \ref{concentration_U.9610} to assert that both $c_n$ and $d_n$ have almost sure limit zero and that in the sums \eqref{12710.13h1} and \eqref{12710.13h2}, any term \st $i\neq j$ tends almost surely to zero.
Moreover, by \eqref{250709.09h13.300709_iid2}, these sums can also be reduced to the terms with index $i$ \st $\tta_i=\tta_{i_0}$. Tu sum up,  we have
 \beq &&a_n=\tta_{i_0}^2\sum_{i\ste \tta_i=\tta_{i_0}}|\lan v_i, {\wtv}\ran|^2 {u_i}^*\f{z_n^2}{(z_n^2 I_n-X_nX_n^*)^2}u_i+o(1)  \\
&&b_n  =\tta_{i_0}^2\sum_{ i\ste \tta_i=\tta_{i_0}}|\lan u_i, {\wtu}\ran|^2{v_i}^*\f{X_n^*X_n}{(z_n^2 I_m-X_n^*X_n)^2}v_i+o(1) \eeq

Now, note that since $z_n$ tends to $\rho$,  \beq \ff{n}\Tr\f{z_n^2}{(z_n^2 I_n-X_nX_n^*)^2}&\ninf &\int\f{\rho^2}{(\rho^2-t^2)^2}\ud \mu_X(t), \\
\ff{m_n}\Tr\f{X_n^*X_n}{(z_n^2 I_m-X_n^*X_n)^2}&\ninf &\int   \f{t^2}{(\rho^2-t^2)^2}\ud \ti{\mu}_X(t),
\eeq hence by Proposition \ref{concentration_U.9610}, almost surely,  \beq a_n&=&\tta_{i_0}^2\int\f{\rho^2}{(\rho^2-t^2)^2}\ud \mu_X(t)\sum_{i\ste \tta_i=\tta_{i_0}}|\lan v_i, {\wtv}\ran|^2+o(1),\\
b_n&  =&\tta_{i_0}^2 \int\f{t^2}{(\rho^2-t^2)^2}\ud \ti{\mu}_X(t)\sum_{ i\ste \tta_i=\tta_{i_0}}|\lan u_i, {\wtu}\ran|^2+o(1).
\eeq
Moreover, by \eqref{250709.09h13.300709_iid3}, for all $i$ \st $\tta_i=\tta_{i_0}$, $$
|\lan u_i, {\wtu}\ran|^2=\tta_{i_0}^2(\vfi_{{\mu_X}_X}(\rho))^2|\lan v_i, {\wtv}\ran|^2+o(1).$$ It follows that $$b_n=\tta_{i_0}^4(\vfi_{{\mu_X}_X}(\rho))^2 \int\f{t^2}{(\rho^2-t^2)^2}\ud \ti{\mu}_X(t)\sum_{i\ste \tta_i=\tta_{i_0}}|\lan v_i, {\wtv}\ran|^2+o(1)$$
Since $a_n+b_n=1+o(1)$, we get $$
\sum_{i\ste \tta_i=\tta_{i_0}}|\lan v_i, {\wtv}\ran|^2\ninf \lf(\tta_{i_0}^2\int\f{\rho^2}{(\rho^2-t^2)^2}\ud \mu_X(t)+\tta_{i_0}^4(\vfi_{{\mu}_X}(\rho))^2 \int\f{t^2}{(\rho^2-t^2)^2}\ud \ti{\mu}_X(t)\ri)^{-1}
$$

The  relations
\beq\tta_{i_0}^2\vfi_{\mu_X}(\rho)\vfi_{\widetilde{\mu}_X}(\rho)&=&1\\
2\int\f{\rho^2}{(\rho^2-t^2)^2}\ud \mu_X(t)&=&\ff{\rho}\vfi_{\mu_X}(\rho)-\vfi_{\mu_X}'(\rho)\\
2\int\f{t^2}{(\rho^2-t^2)^2}\ud \widetilde{\mu}_X(t)&=&-\ff{\rho}\vfi_{ \widetilde{\mu}_X}(\rho)-\vfi_{ \widetilde{\mu}_X}'(\rho)\eeq
 allow to recover the RHS of \eqref{250709.09h13.300709_iid} easily. Via \eqref{250709.09h13.300709_iid3}, one easily deduces \eqref {250709.09h13_iid}.

\section{Proof of Theorems \ref{deloc.sing.vect.250709.rectangular} and \ref{smallest_sing_vect_square}} \label{sec:proof 3} Again, we shall only prove Theorem  \ref{deloc.sing.vect.250709.rectangular}  and suppose the $X_n$'s to be non random.

Let us consider the matrix $M_n(z)$ introduced in Lemma \ref{agl_lin_8710}.
Here, $r=1$, so one easily gets, for each $n$, $$\lim_{z\to+\infty}\det M_n(z)=-\tta^{-2}.$$Moreover, for $b_n:=\si_1(X_n)$ the largest singular value of $X_n$,  looking carefully at the term in $\ff{z^2-b_n^2}$ in $\det M_n(z)$, it appears that  with a \pro which tends to one as $n \longrightarrow \infty $, we have $$\lim_{z\to b_n}\det M_n(z) = +\infty.$$It follows that with a
\pro which tends to one as $n \longrightarrow \infty $, the largest singular value $\wts$ of $\wtX$ is $>b_n$.

Then, one concludes using the second part of Lemma \ref{9710.Bob_Dylan_forever}, as in the proof of Theorem 2.3 of \cite{benaych-rao.09}.

\section{Proof of Theorems \ref{sinfluoutbulk} and \ref{sinfluoutbulksquare2}} \label{sec:proof fluct}
\label{sec.away}
We shall only prove Theorem \ref{sinfluoutbulk}, because Theorem  \ref{sinfluoutbulksquare2} can be proved similarly.

We have supposed that $r=1$. Let us denote $u=u_1$ and $v=v_1$. Then we   have $$P_n=\tta uv^*,$$ with $u\in \K^{n\times 1}$, $v\in \K^{m\times 1}$ random   vectors whose entries are $\nu$-distributed independent random variables, renormalized in the {\it orthonormalized model}, and divided by  respectively $\sqrt{n}$ and $\sqrt{m}$ in the {\it i.i.d.. model}.   We also have  that the matrix  $M_n(z)$  defined in Lemma \ref{agl_lin_8710} is a $2\times 2$ matrix.

Let us   fix an arbitrary $b^*$ \st $b<b^*<\rho$.  Theorem \ref{140709.main.rectangular} implies that almost surely, for $n$ large enough,
$\det [M_n(\cdot)]$ vanishes exactly once in $(b^*, \infty)$. Since moreover, almost surely, for all $n$, $$\lim_{z\to+\infty} \det [M_n(z)]=-\ff{\tta^2}<0,$$ we deduce that almost surely,  for $n$ large enough, $\det [M_n(z)]>0$ for $b^*< z < \widetilde{\si}_1$ and  $\det [M_n(z)]<0$ for $\widetilde{\si}_1<z$.

As a consequence, for  any  real number $x$, for $n$ large enough,
 \beqy\label{300810.12h56} \sqrt{n}(\wts-\rho)< x&\iff& \det  M_n\lf(\rho+\f{x}{\sqrt{n}}\ri)>0.
\eeqy
Therefore,  we have to understand
 the limit distributions of the entries of  $M_n\lf(\rho+\f{x}{\sqrt{n}}\ri)$. They are given by the following
 \begin{lem}For any fixed real number $x$, as $n \longrightarrow \infty $, the distribution of $$\Gamma_n:=\sqrt{n}\lf(M_n\lf(\rho+\f{x}{\sqrt{n}}\ri)-\bbm\vfi_{\mu_X}(\rho)&-\tta^{-1}\\ -\tta^{-1}&\vfi_{\widetilde{\mu}_X}(\rho)\ebm\ri)$$ converges weakly to the one of $$x\bbm \vfi_{\mu_X}(\rho)&0\\ 0&\vfi_{\widetilde{\mu}_X}(\rho)\ebm+\bbm c_1X&dZ\\ d\ovl{Z}&c_2Y\ebm,$$ for $X,Y,Z$ (resp. $X,Y,  \Re(Z),  \Im(Z)$) independent standard real Gaussian   variables   if $\beta=1$ (resp. if $\beta=2$) and for  $c_1,c_2,d$
  some real constants given by the following formulas: \beqy c_1^2&=&\begin{cases} \f{2}{\beta}\int\f{\rho^2}{(\rho^2-t^2)^2}\ud\mu_X(t)&\textrm{in the i.i.d. model,}
  \\
 \f{2}{\beta}\lf(\int\f{\rho^2}{(\rho^2-t^2)^2}\ud\mu_X(t) -  (\vfi_{\mu_X}(\rho))^2\ri)&\textrm{in the orthonormalized model,}
  \end{cases}\label{constant_c_1}\\
   c_2^2&=&\begin{cases} \f{2}{\beta}\int\f{\rho^2}{(\rho^2-t^2)^2}\ud\widetilde{\mu}_X(t)&\textrm{in the i.i.d. model,}
    \\
 \f{2}{\beta}\lf(\int\f{\rho^2}{(\rho^2-t^2)^2}\ud\widetilde{\mu}_X(t) -  (\vfi_{\widetilde{\mu}_X}(\rho))^2\ri)&\textrm{in the orthonormalized model,}\end{cases}\label{constant_c_2}\\
   d^2&=& \ff{\beta}\int\f{t^2}{(\rho^2-t^2)^2}\ud{\mu_X}(t) \label{constant_d}.  \eeqy
 \end{lem}

\begin{pr}
 Let us define $z_n:=\rho+\f{x}{\sqrt{n}}$. We have
$$\Gamma_n= \sqrt{n}\bbm \ff{n}u^*\f{z_n}{z_n^2 I_n-X_nX_n^*}u-\vfi_{\mu_X}(\rho)& \ff{\sqrt{nm_n}}u^*(z_n^2 I_n-X_nX_n^*)^{-1}X_n v \\ \ff{\sqrt{nm_n}}v^*X_n^*(z_n^2 I_n-X_nX_n^*)^{-1}u &\ff{m_n} v^*\f{z_n}{z_n^2 I_m-X_n^*X_n}v -\vfi_{\widetilde{\mu}_X}(\rho)\ebm$$
    Let us for example expand the upper left entry of $\Gamma_{n,1,1}$ of $\Gamma_n$. We have \beqy\nonumber \Gamma_{n,1,1}&=&\sqrt{n}\lf(\ff{n}u^*\f{z_n}{z_n^2 I_n-X_nX_n^*}u-\vfi_{\mu_X}(\rho) \ri)\\
    &=&\sqrt{n}\lf(\ff{n}u^*\f{z_n}{z_n^2 I_n-X_nX_n^*}u-\ff{n}\Tr \f{z_n}{z_n^2 I_n-X_nX_n^*}\ri)\nonumber\\
    &&+\sqrt{n}\lf(\ff{n}\Tr \f{z_n}{z_n^2 I_n-X_nX_n^*}-\vfi_{\mu_X}(z_n)\ri)+\sqrt{n}\lf(\vfi_{\mu_X}(z_n)-\vfi_{\mu_X}(\rho)\ri) \label{29810.13h51} \eeqy
    The third term of the RHS of \eqref{29810.13h51} tends to $x\vfi_{\mu_X}'(\rho)$ as $n \longrightarrow \infty $. By Taylor-Lagrange Formula, there is $\xi_n\in (0,1)$ \st the second one is equal to $$\sqrt{n}\lf(\ff{n}\Tr \f{\rho}{\rho^2 I_n-X_nX_n^*}-\vfi_{\mu_X}(\rho)\ri)+x\f{\partial}{\partial z}_{|_{z=\rho+\xi_n x/\sqrt{n}}}\lf(\ff{n}\Tr \f{z}{z^2 I_n-X_nX_n^*}-\vfi_{\mu_X}(z)\ri),
$$ hence tends to zero, by Assumptions \ref{hypspec} and \ref{H2prime}.
  To sum up, we have  \beqy\Gamma_{n,1,1}&=& \sqrt{n}\lf(\ff{n}u^*\f{z_n}{z_n^2 I_n-X_nX_n^*}u-\ff{n}\Tr \f{z_n}{z_n^2 I_n-X_nX_n^*}\ri)+ x\vfi_{\mu_X}'(\rho)\label{29810.11} +o(1)\eeqy In the same way, we have
  \beqy\Gamma_{n,2,2}&=& \sqrt{n}\lf(\ff{m_n}v^*\f{z_n}{z_n^2 I_m-X_n^*X_n}v-\ff{m_n}\Tr \f{z_n}{z_n^2 I_m-X_n^*X_n}\ri)+ x\vfi_{\widetilde{\mu}_X}'(\rho)\label{29810.22} +o(1)\eeqy  Then the ``$\kappa_4(\nu)=0$" case of Theorem 6.4 of \cite{bg-maida-guionnetTCL} allows to conclude.
\end{pr}

Let us now complete the proof of Theorem \ref{sinfluoutbulk}. By the previous lemma, we have
\beq &
\det   M_n\lf(\rho+\f{x}{\sqrt{n}}\ri)=&\\ & \det \lf( \bbm\vfi_{\mu_X}(\rho)&-\tta^{-1}\\ -\tta^{-1}&\vfi_{\widetilde{\mu}_X}(\rho)\ebm +\ff{\sqrt{n}}\bbm x\vfi_{\mu_X}(\rho)+c_1X_n&dZ_n\\ d\ovl{Z_n}&x\vfi_{\widetilde{\mu}_X}(\rho)+c_2Y_n\ebm \ri) &
\eeq
for some random variables $X_n,Y_n, Z_n$ with converging in distribution to the random variables $X,Y,Z$ of the previous lemma. Using  the relation $\vfi_{\mu_X}(\rho)\vfi_{\widetilde{\mu}_X}(\rho)=\tta^{-2}$, we get
\beq &
 \det   M_n\lf(\rho+\f{x}{\sqrt{n}}\ri)&\\ &= 0+\ff{\sqrt{n}}\lf\{
2x\tta^{-2}+\vfi_{\mu_X}(\rho)c_2Y_n+\vfi_{\widetilde{\mu}_X}(\rho)c_1X_n+\tta^{-1}d(Z_n+\ovl{Z_n})\ri\}+O\lf(\ff{n}\ri)&
\eeq
Thus by \eqref{300810.12h56}, we have
 \beq \lim_{n\to\infty}\Pro\{\sqrt{n}(\wts-\rho)< x\}&=& \lim_{n\to\infty}\Pro\{\det M_n\lf(\rho+\f{x}{\sqrt{n}}\ri)>0\}\\
 &=&\Pro\{   -\f{\tta^2}{2}\lf(\vfi_{\mu_X}(\rho)c_2Y+\vfi_{\widetilde{\mu}_X}(\rho)c_1X+\tta^{-1}d(Z+\ovl{Z})\ri)<x \}.
\eeq
It follows that the distribution of $\sqrt{n}(\wts-\rho)$ converges weakly to the one of $s X$, for $X$ a standard Gaussian random variable on $\R$ and \bes\label{300810.14h35}s^2=\f{\tta^4}{4}\lf( (\vfi_{\widetilde{\mu}_X}(\rho)c_1)^2+
(\vfi_{{\mu_X}_X}(\rho)c_2)^2+4(\tta^{-1}d)^2
\ri) .\ees
One can easily recover the formula given in Theorem \ref{sinfluoutbulk} for $s^2$, using the relation $\vfi_{\mu_X}(\rho)\vfi_{\widetilde{\mu}_X}(\rho)=\tta^{-2}$.

\section{Appendix} \label{sec:appendix}

We now state  the continuity lemma that we use in the proof of Theorem \ref{140709.main.rectangular}. We note that nothing in its hypotheses is random. As hinted earlier, we will invoke it to localize the extreme eigenvalues of $\wtX$.
\begin{lem}\label{roots_continuity_9610}  We suppose the positive real numbers  $\tta_1, \ldots, \tta_r$ to be pairwise distinct.
Let us fix a   real  number $0\le b$ and two analytic functions $\vfi_1, \vfi_2$ defined on $\{z\in \C\ste \Re(z)>0\}\bck[0,b]$ \st for all $i=1,2$, \begin{itemize}
\item[a)] $\vfi_i(z)\in \R\iff z\in \R$,\\
\item[b)] for all $z>b$, $\vfi_i(z)<0$,\\
\item[c)] $\vfi_i(z)\lto0$ as   $|z|\lto \infty$.
\end{itemize}
Let us define the $2r\times 2r$-matrix-valued function$$M(z):=\bbm 0& \Tta^{-1}\\
\Tta^{-1}&0\ebm-\bbm \vfi_1(z)I_r&0\\ 0&\vfi_2(z)I_r\ebm$$ and denote by $z_1>\cdots >z_p$ the $z$'s in $(b,\infty)$ \st $M(z)$ is not invertible, where $p\in \{0, \ldots, r\}$ is the number of $\tta_i$'s \st   $$\lim_{z\downarrow b}\vfi_1(z)\vfi_2(z)>\ff{\tta_i^2}.$$
Let us also consider a sequence
 $0<b_n$
with limit $b$ and, for each $n$, a $2r\times 2r$-matrix-valued function $M_n(\cdot)$, defined on $$\{z\in \C\ste \Re(z)>0\}\bck[0,b_n],$$ which coefficient are analytic functions,  \st
\begin{itemize}
\item[d)] for all $z\notin \R$, $M_n(z)$ is invertible,\\
\item[e)] for all $\eta>0$, $M_n(\cdot)$ converges to the function $M(\cdot)$ uniformly on $\{z\in \C\ste \Re(z)>b+\eta\}$.
\end{itemize}
Then
\begin{itemize}\item[(i)]  there exists $p$ real sequences  $z_{n,1}>\cdots >z_{n,p}$ converging respectively to $z_1,\ldots,z_p$ \st for any $\eps>0$ small enough,  for $n$ large enough, the $z$'s in $(b+\eps, \infty)$ \st $M_n(z)$ is not invertible are  exactly $z_{n,1}, \ldots, z_{n,p}$,\\ \item[(ii)]  for $n$ large enough, for each $i$, $M_n(z_{n,i})$ has rank $2r-1$.\end{itemize}
 \end{lem}
\begin{pr}To prove this lemma, we use the formula  $$\det\bbm xI_r&\Diag(\al_1, \ldots, \al_r)\\
\Diag(\al_1, \ldots, \al_r) &yI_r\ebm=\prod_{i=1}^r(xy-\al_i^2)$$ in the appropriate place
and  proceed as the proof of Lemma 6.1 in \cite{benaych-rao.09}.\end{pr}



We also need the following proposition. The $u_i$'s and the $v_i$'s
  are the random column vectors introduced in Section \ref{def_and_hyp}.
\begin{propo}\label{concentration_U.9610}Let, for each $n$,    $A_n$, $B_n$ be complex $n\times n$, $n\times m$   matrices which operator norms, with respect to the canonical Hermitian structure, are bounded independently of $n$. Then  for any $\eta>0$, there exists $C,\al>0$ \st for all $n$, for all $i,j,k\in \{1, \ldots, r\}$ \st $i\ne j$,
$$\Pro\{|\lan u_i, A_nu_i\ran-\ff{n}\Tr(A_n)|>\eta\textrm{ or }|\lan u_i, A_nu_j\ran|>\eta \textrm{ or }  |\lan u_i, B_nv_k\ran|>\eta\}\le Ce^{-n^\al}.$$
\end{propo}

\begin{pr}In the {\it i.i.d. model}, this result is an obvious consequence of \cite[Prop. 6.2]{bg-maida-guionnetTCL}. In the {\it orthonormalized model}, one also has to use \cite[Prop. 6.2]{bg-maida-guionnetTCL}, which states that the $u_i$'s (the same holds for the $v_i$'s) are obtained from the $n\times r$  matrix $G^{(n)}_u$ with i.i.d. entries distributed according to $\nu$ by the following formula: for all $i=1, \ldots, r$,
$$ u_i=\f{\textrm{$i$th column of $G^{(n)}_u\times (W^{(n)})^T$}}{\|\textrm{$i$th column of $G^{(n)}_u\times (W^{(n)})^T$}\|_2},$$where $W^{(n)}$ is a (random) $r\times r$ matrix \st  for certain positive constants $D, c, \kappa$, for all $\eps>0$ and   all $n$,  $$\Pro\{\|W^{(n)}-I_r\|>\eps\textrm{ or } \max_{1\le i\le r}\lf|\ff{\sqrt{n}}\|\textrm{$i$th column of $G^{(n)}_u\times (W^{(n)})^T$}\|_2-1\ri|>\eps\}\le D(e^{-cn\eps}+e^{-c\sqrt{n}}).$$\end{pr}


\begin{thebibliography}{10}


\bibitem{alter2000singular}
O.~Alter, P.O. Brown, and D.~Botstein.
\newblock {Singular value decomposition for genome-wide expression data
  processing and modeling}.
\newblock {\em Proceedings of the National Academy of Sciences of the United
  States of America}, 97(18):10101, 2000.




\bibitem{anderson84a}
T.W. Anderson.
\newblock {\em An Introduction to Multivariate Statistical Analysis}.
\newblock Wiley, New York, second edition, 1984.

\bibitem{agz09}
G.~Anderson, A.~Guionnet, O.~Zeitouni.
\newblock \emph{An Introduction to Random Matrices}.
\newblock  Cambridge studies in advanced mathematics, {\bf 118} (2009).

\bibitem{agg88}
P.~Arbenz, W.~Gander, and G.~H. Golub.
\newblock Restricted rank modification of the symmetric eigenvalue problem:
  theoretical considerations.
\newblock {\em Linear Algebra Appl.}, 104:75--95, 1988.


\bibitem{bai-silver-95} Z.~D.~Bai, J.~W.~Silverstein On the empirical distriubtion of eigenvalues of a class of
large dimensional random matrices. \emph{J. Multivariate Anal.} (1995), 175--192.


\bibitem{bai-silver-book} Z.~D.~Bai, J.~W.~Silverstein \emph{Spectral analysis of large dimensional random matrices}, Second Edition, Springer, New York, 2009.

\bibitem{bai-yao-TCL}  Z.~D.~Bai, J.-F.~Yao {Central limit theorems for eigenvalues in a spiked population model,} \emph{Ann. I.H.P.-Prob. et Stat.}  Vol. 44 No. 3 (2008), 447--474.



\bibitem{bbp05}
J.~Baik, G.~Ben~Arous, and S.~P{\'e}ch{\'e}.
\newblock Phase transition of the largest eigenvalue for nonnull complex sample
  covariance matrices.
\newblock {\em Ann. Probab.}, 33(5):1643--1697, 2005.



\bibitem{bs06}
J.~Baik, and J.~W.~Silverstein.
\newblock Eigenvalues of large sample covariance matrices of spiked
              population models.
\newblock {\em J. Multivariate Anal.}, 97(6): 1382--1408, 2006.




\bibitem{fbg05.inf.div}
F. Benaych-Georges.
\newblock{Infinitely divisible distributions for rectangular free convolution: classification and matricial interpretation.}
\newblock{{\em Probab. Theory and Related Fields}, 139, 1--2, 2007, 143--189.}





\bibitem{b09}
F.~Benaych-Georges.
\newblock Rectangular random matrices, related convolution.
\newblock \emph{Probab. Theory Related Fields}, 144(3-4):471--515, 2009.

\bibitem{bg07c} F. Benaych-Georges. \emph{Rectangular random matrices, related free entropy and free Fisher's information}.
 \emph{Journal of Operator Theory}. Vol. 23, no. 2 (2010) 447--465.

\bibitem{b09R=C}
F.~Benaych-Georges.
\newblock On a surprising relation between the {M}archenko-{P}astur law,
  rectangular and square free convolutions.
\emph{Ann. Inst. Henri Poincar\'e Probab. Stat.} Vol. 46, no. 3 (2010) 644--652.


\bibitem{bg.sph.int} F. Benaych-Georges.
\newblock{Rectangular $R$-transform at the limit of rectangular spherical integrals},
\emph{J. Theoret. Probab.} Vol. 24, no. 4 (2011) 969--987.


\bibitem{bg-maida-guionnetTCL} F. ~Benaych-Georges, A.~Guionnet and M.~Ma\"\i da. \newblock{Fluctuations of the extreme eigenvalues of finite rank deformations of random matrices}. \emph{Electron. J. Prob.} Vol. 16 (2011), Paper no. 60, 1621--1662.

 \bibitem{bg-maida-guionnetPGD} F. ~Benaych-Georges, A.~Guionnet and M.~Ma\"\i da. \newblock{Large deviations of the extreme eigenvalues of random deformations of matrices}. To appear in \emph{Probab. Theory Related Fields}.

\bibitem{benaych-rao.09}
F.~Benaych-Georges, R.~R.~Nadakuditi.
\newblock  {The eigenvalues and eigenvectors of finite, low rank perturbations of large random matrices},
\newblock   \emph{Adv. in Math.}, Vol. 227, no. 1 (2011) 494--521.



\bibitem{Capitaine2011} M. Capitaine. {Additive/multiplicative free subordination property and limiting eigenvectors of spiked additive deformations of Wigner matrices and spiked sample covariance matrices}. Preprint, 2011.



\bibitem{capcas.iumj} M. Capitaine, M. Casalis. {Asymptotic freeness by generalized moments for Gaussian and Wishart matrices. Application to beta random matrices}.  \emph{Indiana Univ. Math. J.}  53  (2004),  no. 2, 397--431.

\bibitem{capdonmart}  M. Capitaine, C. Donati-Martin.  {Strong asymptotic freeness for Wigner and Wishart matrices}. \emph{Indiana Univ. Math. J.} 56 (2007), no. 2, 767--803.

\bibitem{CDF09} M. Capitaine, C. Donati-Martin, D. F\'eral. {The largest eigenvalues of finite rank deformation of large Wigner matrices: convergence and nonuniversality of the fluctuations}. \newblock\emph{Ann. Probab.} {37} (2009)1--47.

\bibitem{CDF09b} M. Capitaine, C. Donati-Martin, D. F\'eral. {Central limit theorems for eigenvalues of deformations of Wigner matrices}. To appear in   \emph{Annales de L'Institut Henri Poincar\'e}.


\bibitem{CDFF10} M. Capitaine, C. Donati-Martin, D. F\'eral, M. F\'evrier. {Free convolution with a semi-circular convolution and eigenvalues of spiked deformations of Wigner matrices}. To appear in  \emph{Electronic Journal of Probability}.









\bibitem{drineas2004clustering}
P.~Drineas, A.~Frieze, R.~Kannan, S.~Vempala, and V.~Vinay.
\newblock {Clustering large graphs via the singular value decomposition}.
\newblock {\em Machine Learning}, 56(1):9--33, 2004.


\bibitem{edfors2002ofdm}
O.~Edfors, M.~Sandell, J.J. Van~de Beek, S.K. Wilson, and P.O. Borjesson.
\newblock {OFDM channel estimation by singular value decomposition}.
\newblock {\em Communications, IEEE Transactions on}, 46(7):931--939, 2002.

\bibitem{eckart1936approximation}
C.~Eckart and G.~Young.
\newblock {The approximation of one matrix by another of lower rank}.
\newblock {\em Psychometrika}, 1(3):211--218, 1936.


\bibitem{FP07} D. F\'eral, S. P\'ech\'e. {The largest eigenvalue of rank one deformation of large
              {W}igner matrices}.  \emph{Comm. Math. Phys.}
{272} (2007)185--228.

\bibitem{furnas1988information}
G.W. Furnas, S.~Deerwester, S.T. Dumais, T.K. Landauer, R.A. Harshman, L.A.
  Streeter, and K.E. Lochbaum.
\newblock {Information retrieval using a singular value decomposition model of
  latent semantic structure}.
\newblock In {\em Proceedings of the 11th annual international ACM SIGIR
  conference on Research and development in information retrieval}, pages
  465--480. ACM, 1988.

\bibitem{golub1965calculating}
G.~Golub and W.~Kahan.
\newblock {Calculating the singular values and pseudo-inverse of a matrix}.
\newblock {\em Journal of the Society for Industrial and Applied Mathematics:
  Series B, Numerical Analysis}, 2(2):205--224, 1965.



\bibitem{hachemloubatonmnv11} W. Hachem, P. Loubaton, X. Mestre, J. Najim, P. Vallet {A Subspace Estimator for Fixed Rank Perturbations
of Large Random Matrices}. Preprint, 2011.

\bibitem{hp00}
F. Hiai and D. Petz.
\newblock {\em The semicircle law, free random variables and entropy},
  volume~77 of {\em Mathematical Surveys and Monographs}.
\newblock American Mathematical Society, Providence, RI, 2000.

\bibitem{hj85}
R.~A. Horn and C.~R. Johnson.
\newblock {\em Matrix analysis}.
\newblock Cambridge University Press, Cambridge, 1985.

\bibitem{hj91}
R.~A. Horn and C.~R. Johnson.
\newblock {\em Topics in matrix analysis}.
\newblock Cambridge University Press, Cambridge, 1991.

\bibitem{hr07}
D.~C. Hoyle and M.~Rattray.
\newblock Statistical mechanics of learning multiple orthogonal signals:
  asymptotic theory and fluctuation effects.
\newblock {\em Phys. Rev. E (3)}, 75(1):016101, 13, 2007.

\bibitem{jolliffeprincipal}
I.~Jolliffe.
\newblock {Principal component analysis}
\newblock {\em Springer Series in Statistics}, 2002.


\bibitem{johnstone2009statistical}
I.M. Johnstone and D.M. Titterington.
\newblock {Statistical challenges of high-dimensional data}.
\newblock {\em Philosophical Transactions of the Royal Society A: Mathematical,
  Physical and Engineering Sciences}, 367(1906):4237, 2009.

\bibitem{kannan2009spectral}
R.~Kannan and S.~Vempala.
\newblock {\em {Spectral algorithms}}.
\newblock Now Publishers Inc, 2009.


\bibitem{klema2002singular}
V.~Klema and A.~Laub.
\newblock {The singular value decomposition: Its computation and some
  applications}.
\newblock {\em Automatic Control, IEEE Transactions on}, 25(2):164--176, 2002.

\bibitem{jordan1998learning}
M.I. Jordan.
\newblock {\em {Learning in graphical models}}.
\newblock Kluwer Academic Publishers, 1998.

\bibitem{mp67}
V.~A. Mar{\v{c}}enko and L.~A. Pastur.
\newblock Distribution of eigenvalues in certain sets of random matrices.
\newblock {\em Mat. Sb. (N.S.)}, 72 (114):507--536, 1967.

\bibitem{mirsky1960symmetric}
L.~Mirsky.
\newblock {Symmetric gauge functions and unitarily invariant norms}.
\newblock {\em The Quarterly Journal of Mathematics}, vol. 11(1), 50--59, 1960.

\bibitem{Muirhead82}
R.~J. Muirhead.
\newblock {\em Aspects of Multivariate Statistical Theory}.
\newblock Wiley, New York, 1982.



\bibitem{r09}
R.~R. Nadakuditi and J.~W.~Silverstein.
\newblock Fundamental limit of sample generalized eigenvalue based detection of
  signals in noise using relatively few signal-bearing and noise-only samples.
\newblock {\rm IEEE Journal of Selected Topics in Signal Processing}, vol. 4(3), 468--480, 2010.





\bibitem{n08}
B.~Nadler.
\newblock Finite sample approximation results for principal component analysis:
  a matrix perturbation approach.
\newblock {\em Ann. Statist.}, 36(6):2791--2817, 2008.




\bibitem{pan} G. Pan. Strong convergence of the empirical distribution of eigenvalues of
sample covariance matrices with a perturbation matrix \emph{J.  Multivariate Anal.} 101 (2010) 1330--1338.

\bibitem{p07}
D.~Paul.
\newblock Asymptotics of sample eigenstructure for a large dimensional spiked
  covariance model.
\newblock {\em Statist. Sinica}, 17(4):1617--1642, 2007.


\bibitem{sandrinePTRF06}
S.~P\'ech\'e.
\newblock{The largest eigenvalue of small rank
perturbations of Hermitian random matrices}.
\newblock \emph{Probab. Theory Related Fields},
{134}, (2006) 127--173.

\bibitem{rao1964use}
C.R. Rao.
\newblock {The use and interpretation of principal component analysis in
  applied research}.
\newblock {\em Sankhy{\=a}: The Indian Journal of Statistics, Series A},
  26(4):329--358, 1964.

\bibitem{scharf91a}
L.~L. Scharf.
\newblock {\em Statistical Signal Processing: Detection, Estimation, and Time
  Series Analysis}.
\newblock Addison-Wesley, Reading, Massachusetts, 1991.

\bibitem{schmidt86a}
R.~O. Schmidt.
\newblock Multiple emitter location and signal parameter estimation.
\newblock {\em IEEE Trans. on Antennas Prop.}, AP-34:276--280, March 1986.

\bibitem{ss90}
G.~W. Stewart and J.~G. Sun.
\newblock {\em Matrix perturbation theory}.
\newblock Computer Science and Scientific Computing. Academic Press Inc.,
  Boston, MA, 1990.


\bibitem{tipping1999probabilistic}
M.E. Tipping and C.M. Bishop.
\newblock {Probabilistic principal component analysis}.
\newblock {\em Journal of the Royal Statistical Society: Series B (Statistical
  Methodology)}, 61(3):611--622, 1999.

\bibitem{vantrees02a}
H.~L.~Van Trees.
\newblock {\em Detection, Estimation, and Modulation Theory Part IV: Optimum
  Array Processing}.
\newblock John wiley and Sons, Inc., new York, 2002.

\bibitem{wall2003singular}
M.~Wall, A.~Rechtsteiner, and L.~Rocha.
\newblock {Singular value decomposition and principal component analysis}.
\newblock {\em A practical approach to microarray data analysis}, pages
  91--109, 2003.








\end{thebibliography}
\end{document}